\documentclass{article}
\usepackage{amsfonts}

\usepackage{pgf,pgfarrows,pgfnodes,pgfautomata,pgfheaps,pgfshade}
\usepackage[english]{babel}

\begin{document}
\renewcommand{\theequation}{\thesection.\arabic{equation}}
\newcommand{\nor}[1]{\|#1\|}
\newcommand{\R}   {\mathbb{R}}
\newcommand{\N}   {\mathbb{N}}
\newcommand{\C}   {\mathbb{C}}
\newcommand{\Z}   {\mathbb{Z}}
\newcommand{\0}   {\mbox{\bf 0}}
\newcommand{\alert}  {\marginpar{\Large ???}}
\newcommand{\res}    {\mbox{\rm res}}
\newcommand{\eproof} {\framebox{}}
\newcommand{\app }[1]{\stackrel{\mbox{\scriptsize \rm(A.{#1})}}{=}}
\newcommand{\V}     {\mbox{\it V}}
\newcommand{\U}     {\mbox{\it U}}
\newcommand{\vv}     {\mbox{\it v}}
\newcommand{\uu}     {\mbox{\it u}}
\newcommand{\ww}     {\mbox{\it w}}
\renewcommand{\dim}{\textsf{d}}
\newcommand{\Skw}     {\mbox{\rm Skw}\,}
\newcommand{\Sym}    {\mbox{\rm Sym}\,}
\newcommand{\skw}     {\mbox{\rm skw}\,}
\newcommand{\sym}    {\mbox{\rm sym}\,}
\newcommand{\grad}    {\mbox{\rm grad}\,}
\newcommand{\divt}    {\mbox{\rm Div}\,}
\newcommand{\divv}    {\mbox{\rm div}\,}
\newcommand{\Dym}    {\mbox{\scriptsize \rm Dym}}
\newcommand{\ARTICLE}[1]    { \bibitem{#1}}
\newcommand{\BOOK}[1]       { \bibitem{#1}}
\newcommand{\AUTHOR}[1]{{#1}:}
\newcommand{\JOURNAL}[1]{{#1},}
\newcommand{\TITLE}[1]{{\it{#1}},}
\newcommand{\VOLUME}[1]{{#1},}
\newcommand{\PAGES}[1]{p. {#1},}
\newcommand{\PUBLISHER}[1]{{#1},}
\newcommand{\ADdsESS}[1]{{#1},}
\newcommand{\YEAR}[1]{{#1}.}

\newcounter{mycounter}[section]
\def\themycounter{\thesection.\arabic{mycounter}}
\newcommand{\definition}[1]{\refstepcounter{mycounter}      \vspace{3.0ex}
   \noindent{\bf Definition \themycounter~~} {\sl #1}\par   \vspace{3.0ex}}
\newcommand{\lemma}[1] {   \refstepcounter{mycounter} % ripristinato da G.Floridia
 \vspace{3.0ex}
   \noindent{\bf Lemma \themycounter~~} {\sl #1}\par     \vspace{3.0ex}}
\newcommand{\theorem}[1]{\refstepcounter{mycounter}      \vspace{3.0ex}
   \noindent{\bf Theorem \themycounter~~} {\sl #1}\par   \vspace{3.0ex}}
\newcommand{\remark}[1]{\refstepcounter{mycounter}       \vspace{3.0ex}
   \noindent{\bf Remark }%\themycounter~~
 {\rm #1 {~}\hfill \eproof\par\goodbreak\vspace{3.0ex}}}
%}\par    \vspace{3.0ex}}
\newcommand{\proposition}[1]{\refstepcounter{mycounter}       \vspace{3.0ex}
   \noindent{\bf Proposition \themycounter~~} {\sl #1}\par    \vspace{3.0ex}}
\newcommand{\example}[1]{\refstepcounter{mycounter}      \vspace{3.0ex}
   \noindent{\bf Example \themycounter:~} {#1}\par   \vspace{3.0ex}}
\newcommand{\Problem}[1]{\refstepcounter{mycounter}      \vspace{3.0ex}
   \noindent{\bf Problem \themycounter :~} {#1}\par   \vspace{3.0ex}}
\newcommand{\proof}[1]{\vspace{-0.1cm}\noindent {\bf Proof: {\rm #1 }}  \newline }
\bigskip\bigskip\bigskip
\begin{center} {\Large \bf
{ Approximate controllability for linear degenerate parabolic problems with bilinear
control}
%{ Title in one line if short enough or, if not, \\ in more lines}
 }
\end{center}
\normalsize \vspace{0.2cm}

\medskip
\begin{center}
%{\bf Author 1,
%\\{\rm Dipartimento ,\\
%        Universit\`a di Roma ``La Sapienza'',\\
%        I-00161 Roma, Italy}
% \\ \medskip
%{\bf Author 2,
%\\{\rm Dipartimento ,\\
%        Universit\`a di Roma ``La Sapienza'',\\
%        I-00161 Roma, Italy}
% }}
{\bf Piermarco Cannarsa,
\\{\rm Dipartimento di Matematica ,\\
        Universit\`a di Roma ``Tor Vergata'',\\
        I-00161 Roma, Italy}
 \\ \medskip
{\bf Giuseppe Floridia,
\\{\rm Dipartimento di Matematica e Informatica ,\\
        Universit\`a di Catania,\\
        I-95125 Catania, Italy}
 }}
\end{center}
\medskip
%%%%%%%%%%%%%%%%%%%%%%%%%%%%%%%%%%%%%%%%%%%%%%%%%%%%%%%%%%%%%%%%%%%%%%%%%%%%%%%%%%%%%%%%%%%%%%%%%%%%%%
\begin{abstract}
In this work we study the global approximate multiplicative controllability for the linear degenerate parabolic Cauchy-Neumann problem
$$
\left\{\begin{array}{l}
\displaystyle{v_t-(a(x) v_x)_x =\alpha (t,x)v\,\,\qquad \mbox{in} \qquad Q_T \,=\,(0,T)\times(-1,1) }\\ [2.5ex]
\displaystyle{a(x)v_x(t,x)|_{x=\pm 1} = 0\,\,\qquad\qquad\qquad\,\, t\in (0,T) }\\ [2.5ex]
\displaystyle{v(0,x)=v_0 (x) \,\qquad\qquad\qquad\qquad\quad\,\, x\in (-1,1)}~,
\end{array}\right.
$$
with the bilinear control $\alpha(t,x)\in L^\infty (Q_T).$ The problem is strongly degenerate in the sense that $a\in C^1([-1,1]),$ positive on $(-1,1),$ is allowed to vanish at $\pm 1$ provided that a certain integrability
condition is fulfilled.
We will show that the above system can be steered in $L^2(\Omega)$ from any nonzero, nonnegative initial state into any neighborhood of any desirable nonnegative target-state by bilinear static controls. % (x-dependent only).
%In the particular case $a(x)=1-x^2,$ the above system represents the Budyko-Sellers one-dimensional climatology model.
%In this case we study the controllability properties in large time relating to understand man's possible actions on the environment, in order to intervene the evolution of the temperature.
Moreover, we extend the above result relaxing the %nonnegative
sign constraint on $v_0$.
\end{abstract}
%%%%%%%%%%%%%%%%%%%%%%%%%%%%%%%%%%%%%%%%%%%%%%%%%%%%%%%%%%%%%%%%%%%%%%%%%%%%%%%%%%%%%%%%%%%%%%%%%%%%%%%%%%%%%%%%%%


%%%%%%%%%%%%%%%%%%%%%%%%%%%%%%%%%%%%%%%%%%%%%%%%%%%%%%%%%%%%%%%%%%%%%%%%%%%%%%%%%%%%%%%%%%%%%%%%%%%%%%%%%%%%%%%%%%

\section{Introduction}
\subsection*{Motivation}
Climate depends on various parameters such as temperature,
humidity, wind intensity, the effect of greenhouse gases, and so
on. It is also
%and is
affected by a complex set of interactions in the atmosphere,
oceans and continents, %which they speak of
 that involve  physical, chemical,
geological and biological processes. %\\
%% The weather is also involved in
% Climatology deals with the problem to understand
%% understanding the problem of
% global warming on the Planet.

 One of the first attempts to model the effects of interaction between large ice masses and solar radiation on
climate is the one due, independently, by Budyko \cite{B1, B2}
and Sellers \cite{S} %, the mathematical analysis of those
%models can be shown, for example, in
 (see also \cite{D, DH, H} and the
references therein).  Such a model studies
how extensive the climate response is to an event such as a sharp
increase in greenhouse gases; in this case we talk about climate
sensitivity. A process that changes climate sensitivity is called
\textit{feedback}. If the process increases the intensity of
response we say that it has
 \textit{positive feedback}, whereas it has  \textit{negative feedback} if it reduces the intensity of response.

 The Budyko-Sellers  model studies the role played by
continental and oceanic areas of ice on climate change.
%At the
%beginning of an ice age, the climate cools, then the ice on the
%continents and oceans gains. The average albedo of the Earth and
%significantly increases the amount of energy absorbed by the
%surface.\\
In such a model, the sea level mean
zonally averaged temperature $u(t, x)$ on the Earth, where $t$
denotes time and $x$ the sine of latitude, satisfies the
following degenerate \textit{Cauchy-Neumann} problem (\ref{BS}) in
the bounded domain $(-1,1)$.
%\vspace{-0.5cm}
%\begin{equation}\label{1}
%u_t-(a(x) u_x)_x =\alpha (t,x)u + f (t,x,u)\,\,\qquad \mbox{in}
%\qquad Q_T \,:=\,(0,T)\times(-1,1),
%\end{equation}
%\vspace{-0.5cm}
%\begin{equation}\label{2}
%a(x)u_x(t,x)|_{x=\pm 1} = 0, \,\,\qquad\qquad\quad t\in (0,T),
%\end{equation}
%\vspace{-0.5cm}
%\begin{equation}\label{3}
%\,u(0,x)=u_0 (x), \,\,\qquad\qquad\qquad\qquad x\in (-1,1),
%\end{equation}
%%%%%%%%%%%%%%%%%%%%%%%%%%%%%%%%%%%%%%%%%%%%%%%%%%%%%%%%%%%%%%%%%%

The effect of solar radiation on climate can be summarized in the following figure

\begin{figure}[htbp]
\begin{center}
\hskip-.5cm
\includegraphics[width=5.5cm, height=6cm]{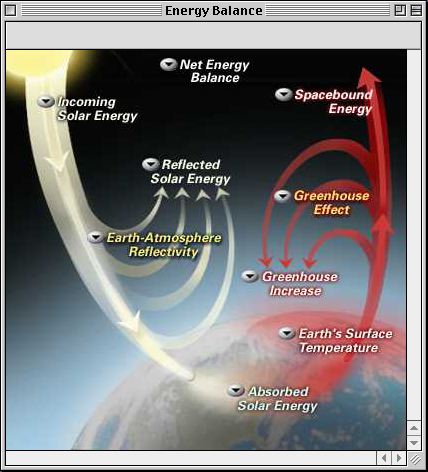}
\caption{www.edu-design-principles.org (copyrighted by DPD)}
%\label{default}
\end{center}
\end{figure}
We have the following {\it energy balance }:\hspace{0.5cm}
\fbox{
%\rule[-.5cm]{0cm}{1cm}
\parbox{1.6in}{
$$\!\!\!\!\mbox{Heat variation}=R_a-R_e+D$$
{\small\begin{itemize}
\item $ R_a$ = absorbed
energy
\item $ R_e$ = emitted energy
\item $ D$ = diffusion
\end{itemize}}
}}
%%%%%%%%%%%%%%%%%%%%%%%%%%%%%%%%%%%%%

%%%%%%%%%%%%%%%%%%%%%%%%%%%%%%%%%%%%%%%%%%%%%%%%%%%%%%%%%%%%%%%%%%
The general formulation of the Budyko-Sellers model on a compact surface $\mathcal M$ without boundary is as follows\\
%considering
 %compact surface without boundary (typically $ S^2$),\\
%then the following equation is obtained
$$u_t-\Delta_{\mathcal M} u= R_a(t,x,u)-R_e(t,x,u)$$

\noindent where\quad $u(t,x)$ is the distribution of temperature and $\Delta_{\mathcal M}$ is the classical Laplace-Beltrami operator. %and
Moreover, %we have
%\pause
\begin{itemize}
\item $ R_a(t,x,u)=Q(t,x)\beta(x,u)$
\item $ R_e(t,x,u)=A(t,x)+B(t,x)u$
\end{itemize}
In the above,
$ \,Q$ is the \textit{insolation} function,
%\qquad\qquad
and
$\,\beta$ is the \textit{coalbedo} function (that is, 1-\textit{albedo} function).
\quad
%\pause
\\
%%%%%%%%%%%%%%%%%%%%%%Wikipedia
%Albedo, or reflection coefficient, is the diffuse reflectivity  or reflecting power of a surface. It is defined as the %ratio of reflected radiation from the surface to incident radiation upon it. Being a dimensionless  fraction, it may %also be expressed as a percentage, and is measured on a scale from zero for no reflecting power of a perfectly black %surface, to 1 for perfect reflection of a white surface.\\
%%%%%%%%%%%%%%%%%%%%%%%%%%%%%%%%%%%%%%%%%%
Albedo is the reflecting power of a surface. It is defined as the ratio of reflected radiation from the surface to incident radiation upon it. It may also be expressed as a percentage, and is measured on a scale from zero for no reflecting power of a perfectly black surface, to 1 for perfect reflection of a white surface.\\

\begin{figure}[htbp]
\begin{center}
%\hskip-.5cm
\includegraphics[width=4cm, height=5cm]{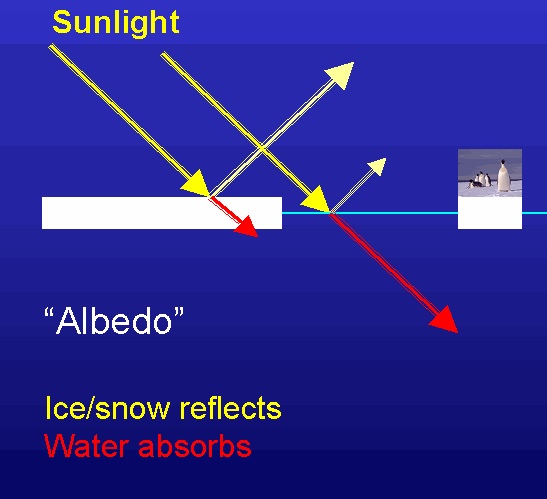}
\caption{www.esr.org (copyrighted by ESR)}
%\label{default}
\end{center}
\end{figure}
%%%%%%%%%%%%%%%%%%%%%

The main difference between Budyko's model and the one by Sellers, %model,
is that in
the former %Budyko
the coalbedo function is discontinuous, while %that of Sellers
in the latter it is a
continuous function. In fact we have\\

%%%%%%%%%%%%%%%
\fbox{
%\rule[-.5cm]{0cm}{1cm}
\parbox{2in}{
\begin{itemize}
\item Budyko
$$
\hspace{-.5cm}\beta(u)=
\left\{\begin{array}{l}
\displaystyle{\beta_0
\,\,\qquad u<-10 }\\ [2.5ex]
\displaystyle{[\beta_0,\beta_1]\quad u=-10  }\\ [2.5ex]
\displaystyle{\beta_1 \,\,\qquad u>-10 }~,
\end{array}\right.
$$
%\pause
\item Sellers
$$
\hspace{-.5cm}\beta(u)=
\left\{\begin{array}{l}
\displaystyle{\beta_0
\,\,\qquad u<u_- }\\ [2.5ex]
\displaystyle{\mbox{line}\quad u_-\leq u\leq u_+  }\\ [2.5ex]
\displaystyle{\beta_1 \,\,\qquad u>u_+ }~,
\end{array}\right.
$$
where $u_\pm=-10\pm\delta, \delta>0.$
\end{itemize}
}}
%%%%%%%%%%%%%%%

On \quad $\mathcal M=\Sigma^2$ the Laplace-Beltrami operator is
$$\Delta_{\mathcal M}=\frac1{\sin \phi}\Big\{\frac{\partial}{\partial \phi}\Big(\sin \phi \frac{\partial u}{\partial\phi}\Big)+\frac1{\sin \phi}\,\frac{\partial^2u}{\partial \lambda^2}\Big\}$$
where $\phi$ is the \textit{colatitude} and $\lambda$ is the \textit{longitude}.\\
%%%%%%%%%%%%%%%%%%%%%%%%%%%%%%%%

\begin{figure}[htbp]
\begin{center}
%\hskip-.5cm
\includegraphics[width=4cm, height=3cm]{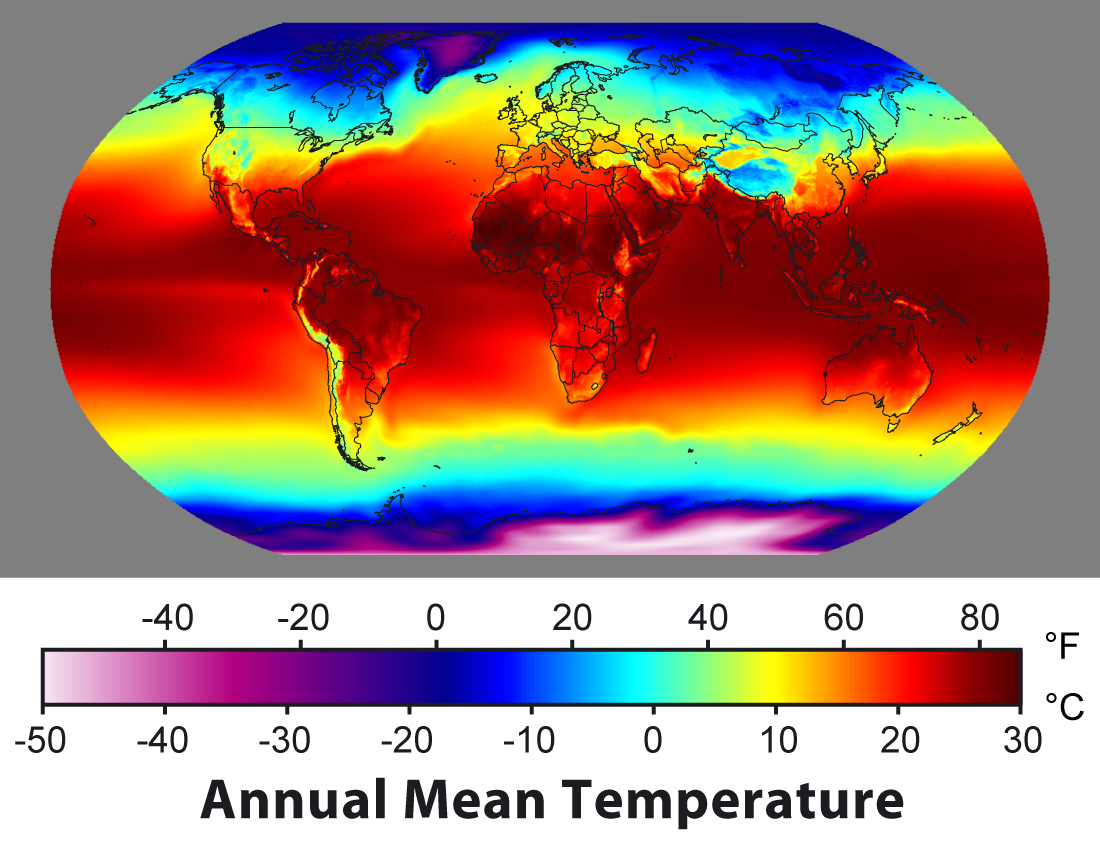}
\caption{www.globalwarmingart.com (copyrighted by Global Warming Art)}
%\label{default}
\end{center}
\end{figure}

%%%%%%%%%%%%%%%%%%%%%%%%%%%%%%%%
%%%%%%%%%%%%%%%%%%%%%%%%%%%%%%%%%
%\begin{figure}[htbp]
%\begin{center}
%%\hskip-.5cm
%\includegraphics[width=4cm, height=3cm]{Annual_Average_Temperature_Map}
%%\caption{default}
%%\label{default}
%\end{center}
%\end{figure}
%%\pause
%\vspace{-.3cm}
%%%%%%%%%%%%%%%%%%%%%%%%%%%%%%%%%%
\vspace{0.5cm}
In the one-dimensional Budyko-Sellers we
take the average of the temperature at $ x=\cos\phi$ and the Budyko-Sellers model reduces to
%$$u_t-\big((1-x^2)u_x\big)_x=g(t,x)\,h(x,u)+f(t,x) \qquad x\in ]-1,1[$$
%\vspace{.5mm}
%\\
%$$(1-x^2)u_{x|_{x=\pm 1}}=0$$
%\end{cases}
%\end{equation*}

\begin{equation}\label{BS}
  \left\{\begin{array}{l}
\displaystyle{u_t-\big((1-x^2)u_x\big)_x=g(t,x)\,h(x,u)+f(t,x), \qquad x\in (-1,1)}\\ [2.5ex]
\displaystyle{(1-x^2)u_{x|_{x=\pm 1}}=0 }~.
\end{array}\right.
\end{equation}

%%%%%%%%%%%%%%%%%%%%%%%%%%%%%%%%%%%%%%%%%%%%%%%%%%%%%%%%%%%%%%%%%%%%%%%%%%%%%%%%%%%%%%%%%%%%%%%%%%%%%%%%%%%%%%%

\subsection*{Problem formulation}

Let us consider the following \textit{Cauchy-Neumann}
strongly degenerate  boundary linear problem in divergence form,
 governed in the bounded domain $(-1,1)$ by means of the \textit{bilinear control} $\alpha (t,x)$ %\in L^\infty (Q_T)$
\\
%\vspace{1.5cm}
\begin{equation}
\label{P2}
%{\displaystyle{
%	\cases{\displaystyle\  v_t-(a(x) v_x)_x =\alpha (t,x)v  ~~~~~ \\\\  &$(t,x)\in Q_T
%%\,=\,(0,T)\times(-1,1),
% $\cr \\
%\displaystyle\ a(x)v_x(t,x) %\stackrel{x \rightarrow\pm 1}{\longrightarrow}
%|_{x=\pm 1} = 0  ~~~~~ \\\\ &$t\in (0,T) $\cr \\
%\displaystyle\ v(0,x)=v_0 (x)  ~~~~~ \\\\  &$x\in (-1,1), $\cr
%}}}
%%%%%%%%%%%%%%%%%%%
\left\{\begin{array}{l}
\displaystyle{v_t-(a(x) v_x)_x =\alpha (t,x)v\,\,\qquad \mbox{in} \qquad Q_T \,=\,(0,T)\times(-1,1) }\\ [2.5ex]
\displaystyle{a(x)v_x(t,x)|_{x=\pm 1} = 0\,\,\qquad\qquad\qquad\,\, t\in (0,T) }\\ [2.5ex]
\displaystyle{v(0,x)=v_0 (x) \,\qquad\qquad\qquad\qquad\quad\,\, x\in (-1,1)}~.
\end{array}\right.
%%%%%%%%%%%%%%%%
%\begin{cases}
%v_t-(a(x) v_x)_x =\alpha (t,x)v,\,\,\qquad \mbox{in} \qquad Q_t \,=\,(0,T)\times(-1,1),
%\\
%a(x)v_x(t,x)
%%\stackrel{x \rightarrow\pm 1}{\longrightarrow}
%|_{x=\pm 1} = 0, \,\,\qquad\qquad\qquad\,\, t\in (0,T),
%\\
%v(0,x)=v_0 (x) \,\qquad\qquad\qquad\qquad\qquad\, x\in (-1,1),
%\end{cases}
\end{equation}

%%%%%%%%%%%%%%%%%%%%%%%%%
\noindent We assume that
\begin{enumerate}
  \item $v_0 \in L^2(-1,1)$
  \item $\alpha \in L^\infty (Q_T)$
  \item $a \in C^1([-1,1])$ satisfies
  \begin{enumerate}
    \item $a(x)>0 \,\, \forall \, x \in (-1,1),\quad a(-1)=a(1)=0$
    \item $A\in L^1(-1,1),$ where $A(x)=\int_0^x \frac{ds}{a(s)}\,.$
  \end{enumerate}
        \end{enumerate}
%where $v_0 (x)\in L^2(-1,1)$, the bilinear coefficient $\alpha (t,x)\in L^\infty (Q_T)$ and the nonnegative function %$a(x):[-1,1]\longrightarrow \R$ that satisfies the following properties
%$$a \in C^1([-1,1]), \quad a(x)>0 \,\, \forall \, x \in
% (-1,1),\quad a(1)=a(-1)=0\,\,(\footnote{ We observe that $\frac{1}{a}\not\in L^1(-1,1),$ so $a(x)$ is strongly %degenerate.}).$$

{\remark{
We observe that
\begin{enumerate}
  \item %We observe that
  $\frac{1}{a}\not\in L^1(-1,1),$ so $a(x)$ is strongly degenerate
  \item %We observe that
   the principal part of the operator in (\ref{P2}) coincides with that of the Budyko-Sellers model for
$a(x)=1-x^2$. In this case $A(x)=\frac{1}{2}\ln\left(\frac{1+x}{1-x}\right)\in L^1(-1,1)$
\item a sufficient condition for 3.b) is that $a^\prime(\pm1)\neq0$ (if $a\in C^2([-1,1])$ the above condition is also necessary).
%\begin{itemize}
%  \item if $a\in C^2([-1,1])$ the above condition is also necessary.
%\end{itemize}
\end{enumerate}
}}
% We observe that the principal part of the operator in (\ref{P2}) coincides with that of the Budyko-Sellers model for
%$a(x)=1-x^2$.
 We are interested in studying the multiplicative controllability of
problem (\ref{P2}) by the \textit{bilinear control} $\alpha (t,x)%\in L^\infty (Q_T)
$. In particular, for the
above linear problem,  we will discuss results guaranteeing
global nonnegative approximate controllability in large time (for
multiplicative controllability see \cite{K,KB,CK}).\\
%\vspace{0.5cm}
Now we recall %the general
one definition %as \textit{globally approximately controllable system.}
from control theory.
%\vspace{1.5cm}

%{\definition{\label{1.1}}{ \it
%We say that the system %\eqref{P1} or
%(\ref{P2}) is globally approximately controllable in the given linear
%phase-space $H$ at time $T>0$, if it can be steered from any
%initial state in $H$ into any neighborhood of any desirable target
%state in $H$ at time $T$, by selecting a suitable available
%control.}}

%\vspace{0.5cm}
%Now we give the particular  definition as \textit{nonnegatively} globally approximately controllable system.

{\definition{\label{1.2}}{ \it
We say that the system (\ref{P2}) %(alternately \eqref{P2})
is \textit{nonnegatively} globally approximately
controllable in $L^2 (-1,1),$ if for every $\varepsilon >0$ and
for every \textit{nonnegative} %$u_0(x),\, u_d(x)\in L^2(\Omega)$ with $u_0\neq 0 \,(
$v_0 (x),\,v_d(x)\in L^2(-1,1)\mbox{ with }
v_0 \not\equiv 0 %)
$ there are a
$%T=T(\varepsilon,u_0,u_d)\,(
T=T(\varepsilon,v_0,v_d)%)
$ and a
bilinear control $\alpha(t,x)\in L^\infty (Q_T)$ such that for the
corresponding solution $%u(t,x)\,(
v(t,x)%)
$ of %\eqref{P1}
%(alternately
(\ref{P2}) %)
 we obtain
%$$\|u(T,\cdot)-u_d\|_{L^2(\Omega)}\leq \varepsilon$$
$$%\mbox{or  }
\|v(T,\cdot)-v_d\|_{L^2(-1,1)}\leq \varepsilon\,.$$}}

%%%%%%%%%%%%%%%%%%%%

In the following, we will sometimes use $\|\cdot\|$ %$\Omega$
 instead of
%the bounded open interval (-1,1)
$\|\cdot\|_{L^2(-1,1)}$. \\

%%%%%%%%%%%%%%%%%%%%%%%%%%%%%%%%%%%%%%%%%%%%%%%%%%%%%%%%%%%%%%%%%%%%%%%%%%%%%%%%%

\subsection*{Main results}
In this work at first the \textit{nonnegative global approximate
controllability} result is obtained for the linear system
(\ref{P2}) in the following theorem. %\ref{T1}
%and then, using this
%result we prove the nonnegatively globally approximately
%controllability for the semilinear problem \ref{P1} in
%\thmref{T2}.

{\theorem{\label{T1}}{
\it The linear system (\ref{P2}) is nonnegatively
approximately controllable in $L^2(-1,1)$ %in the sense of Definition \ref{1.2}
by means of static controls in %$\alpha=\alpha(x), \,
%\alpha\in
$L^\infty (-1,1)$. Moreover, the corresponding
solution to (\ref{P2}) remains nonnegative at all times.
}}

%\begin{thm}
%\label{T2} The semilinear system \eqref{P1} is nonnegatively
%approximately controllable in $L^2(-1,1)$ in the sense of
%\defnref{1.2}. The target state can be achieved by subsequent applying
%of three suitable static bilinear controls.
%\end{thm}

%I risultati presenti nel \thmref{1.1} e \thmref{1.2} possono
%essere estesi a un più grande insieme di stati iniziali, cioè non
%necessariamente non negativi e precisamente alle coppie stato
%iniziale/stato finale $(u_0,u_d)\, or \,(v_0,v_d)$ tali che
\vspace{0.5cm}
Then the results present in Theorem \ref{T1} %and \thmref{1.2}
can be
extended to a larger class of initial states. %, therefore not
%necessarily nonnegative and which are precisely the initial state
%$v_0$ such that
%\begin{equation}
%\label{H2}
%\int^1_{-1}v_0 v_d dx>0
%\end{equation}
%where $v_d%\geq 0,
%$ is any nonnegative target state.
%%in the case of problem \eqref{P2},
%%\\
%%and
%%\begin{equation} \label{H1}
%%\int^1_{-1}u_0 u_d dx>0
%%\end{equation}

%{\theorem{\label{C1}}{ \it
%%The linear system (\ref{P2}) is globally approximately
%%controllable, in the sense of Definition \ref{1.1}, in the space $H$ of
%%initial states $v_0 \in L^2 (-1,1)$ such that (\ref{H2}) holds,
%%with $v_d\in L^2 (-1,1),\, v_d\geq 0$ in $\Omega$.
%%%%%%%%%%%%%%%%%%%%%%%
%The linear system (\ref{P2}) is globally approximately
%controllable in the way that can be steered in $L^2 (-1,1)$ from any initial states $v_0 \in L^2 (-1,1)$ such that %(\ref{H2}) holds into any neighborhood of any desirable nonnegative target state $v_d\in L^2 (-1,1).$}}
%%$,\, v_d\geq 0$ in $\Omega$.

{\theorem{\label{C1}}{ \it
\noindent For any $%\forall
 v_d\in L^2 (-1,1), v_d\geq 0$ and any $%\forall
 v_0\in L^2 (-1,1)$ such that
\begin{equation}
\label{H2}
\int^1_{-1}v_0 v_d dx>0,
\end{equation}
for every $%\forall
\varepsilon>0,$ there are %\exists\,
$ T=T(\varepsilon,v_0,v_d)\geq 0$
and a static bilinear control, %\exists
$\,\alpha=\alpha (x),\,\alpha\in L^\infty(-1,1)$ such that
$$\|v(T,\cdot)-v_d\|_{L^2(-1,1)}\leq \varepsilon\,.$$
%The linear system (\ref{P2}) is globally approximately
%controllable in the way that can be steered in $L^2 (-1,1)$ from any initial states $v_0 \in L^2 (-1,1)$ such that %(\ref{H2}) holds into any neighborhood of any desirable nonnegative target state $v_d\in L^2 (-1,1).$
}}

{\remark{
%La soluzione $v(t,x)$ del problema \eqref{P2}, nelle ipotesi del
%\corref{C1} non rimane non negativa in $Q_T$, come nel
%\thmref{T1}, ma può assumere anche valori negativi.
The solution $v(t,x)$ of the problem (\ref{P2}) in the
assumptions of Theorem \ref{C1} does not remain nonnegative in $Q_T$,
like in Theorem \ref{T1}, but it can also assume negative values.
}}

%\begin{cor}
%\label{C2} The semilinear system \eqref{P1} is globally
%approximately controllable, in the sense of \defnref{1.1}, in the
%space $H$ of initial states $u_0 \in L^2 (-1,1)$ such that
%\eqref{H1} holds; with $u_d\in L^2 (-1,1),\, u_d\geq 0$ in
%$\Omega$.
%\end{cor}

%%%%%%%%%%%%%%%%%%%%%%%%%%%%%%%%%%%%%%%%%%%%%%%%%%%%%%%%%%%%%%%%%%%%%%%%%%%%%%%%%%%%%%%%%%%%%%%%%%%%%%%%%%%%%%%

\subsection*{Mathematical motivation}
This note is inspired by \cite{K,CK}. %in particular
In \cite{K} A.Y. Khapalov studied the global nonnegative approximate controllability of the one dimensional \textit{non-degenerate} semilinear convection-diffusion-reaction equation governed in a bounded domain via the bilinear control $\alpha\in L^\infty (Q_T).$ %in the additive reaction term.
In \cite{CK}, the same approximate controllability property is derived in suitable classes of functions that change sign.\\
%In the linear case Khapalov studies the following equation
%$v_t=v_{xx}+\alpha v%-f(t,x,u,u_x)
%,\,\, \mbox{ in } Q_T=(0,T)\times(-1,1)$ with Dirichlet boundary condition $v(t,0)=v(t,1),\,\,t\in (0,T).$
%The author obtains the nonnegativity of the solution of the previous problem through the use the \textit{maximum %principle} (see also \cite{LSU}). In contrast the use of \textit{maximum principle} is not possible in the degenerate %case with Neumann boundary condition. However in this work this difficulty is surmounted by using an different %technique which we find in Lemma \ref{NN}. This technique is a general method to prove the nonnegativity the solution %v(t,x) in $Q_T$ (in fact we prove by energy estimate and Gronwall's inequality that $v^-$, the negative part of the %function $v$, is equals to zero in $Q_T$ ), that could also be used in the non-degenerate case. \\
In this note we extend some of the results of %multiplicative controllability obtained by Khapalov from non-degenerate case
\cite{K} to degenerate linear equations.\\ %For
General references for \textit{multiplicative controllability} %see also the following references
are, e.g., \cite{K1,K2,K3,K4,KB,BS}.\\
In control theory, boundary and interior locally distributed controls are usually employed (see, e.g., \cite{CMV2,CMV1,CMV3,FC,FCZ,FI}).
These controls are additive terms in the equation and have localized support. %As regards of applications it seem that these controls cannot model all physical process but only process that do not change their principal physical characteristics due to the control action.
%In this way many new technologies are excluded,
However, such models are unfit to study several interesting applied problems such as chemical reactions controlled by catalysts, and also smart materials, which are able to change their principal parameters under certain conditions. %From this the growing interest in the study of the \textit{multiplicative controllability}.
This explains the growing interest in \textit{multiplicative controllability}.

%%%%%%%%%%%%%%%%%%%%%%%%%%%%%%%%%%%%%%%%%%%%%%%%%%%%%%%%%%%%%%%%%%%%%%%%%%%%%%%%%

\section{Preliminaries}
\label{cap1}
\subsection*{Positive and negative part}\label{parti}
Given $\Omega\subseteq\R^n$, $v:\Omega\longrightarrow\R$ we consider the positive-part function
$$
v^+(x) = \max\left(v(x),0\right),\qquad\qquad\qquad\forall x\in \Omega\,,
$$
and the negative-part function
$$
\!\,v^-(x) = \max\left(0,-v(x)\right),\qquad\qquad\qquad\forall x\in \Omega\,.%(\footnote{Calls respectively positive and negative part to the function $v$.})
$$
Then we have the following equality
$$
v=v^+ -v^- \qquad\quad\quad \mbox{  in    }\,\Omega
$$
For the functions $v^+$ and $v^-$ the following result of
regularity in Sobolev's spaces will be useful (see \cite{KS},
Appendix $A$ %, and \cite{M}, Chapter 3
).\\

\vspace{-0.75cm}

{\theorem{\label{A.1}}{ \it
Let $\Omega\subset\R^n, \, u:\Omega\longrightarrow\R,\,  u\in H^{1,s}(\Omega),
\,1\leq s\leq\infty$. Then
$$%\max\,\{u,0\}
u^+,\,u^-\in H^{1,s}(\Omega)$$
and for $1\leq i\leq n$
\begin{equation}
%[\,\max(u,0)\,]
(u^+)_{x_i}=
\left\{\begin{array}{l}
{ u_{x_i} \qquad\qquad\qquad\, \mbox{ in }\{x\in\Omega : u(x)>0\}}\\ [2.5ex]
{ 0 \qquad\qquad\qquad\,\,\,\,\mbox{ in }\{x\in\Omega : u(x)\leq0\} }~,
\end{array}\right.
%(u^+)_{x_i}=
%\begin{cases}
%u_{x_i} \qquad\qquad\qquad\, \mbox{ in }\{x\in\Omega : u(x)>0\}
%\\
%0 \qquad\qquad\qquad\quad \mbox{ in }\{x\in\Omega : u(x)\leq0\}
%\end{cases}
\end{equation}
and
\begin{equation}
%%[\,\max(u,0)\,]
(u^-)_{x_i}=\left\{\begin{array}{l}
{ -u_{x_i} \qquad\qquad\qquad\!\! \mbox{ in }\{x\in\Omega : u(x)<0\}}\\ [2.5ex]
{ 0 \,\,\qquad\qquad\qquad\,\,\, \mbox{ in }\{x\in\Omega : u(x)\geq0\} }~.
\end{array}\right.
\end{equation}
%\begin{cases}
%-u_{x_i} \qquad\qquad\qquad\, \mbox{ in }\{x\in\Omega : u(x)<0\}
%\\
%0 \qquad\qquad\qquad\qquad\, \mbox{ in }\{x\in\Omega : u(x)\geq0\}
%\end{cases}
}}

\subsection*{Gronwall's Lemma}
{\lemma{Gronwall's inequality (differential form).}{\label{GW}}{ \it
Let $\eta(t)$ be a nonnegative, absolutely continuous
function on $[0,T]$, which satisfies for a.e. $t\in [0,T]$ the differential
inequality
\begin{equation}\label{GWI}
\eta'(t)\leq \phi (t) \eta(t)+\psi(t),
\end{equation}
where $\phi (t)$ and $\psi(t)$ are nonnegative, summable functions
on $[0,T]$.
\\
Then
\begin{equation}\label{GWT} \eta(t)\leq e^{\int^t_0 \phi
(s)ds}\left[\eta(0)+\int^t_0 \psi (s)ds\right]
\end{equation}

for all $0\leq t\leq T$.
\\
In particular, if $\psi(t)\equiv 0$ in (\ref{GWI}), i.e.
$
\eta^\prime\leq\phi\,\eta \mbox{ for a.e. } t\in [0,T],
$
and
$\eta(0)=0,$
then
$$
\eta\equiv 0\qquad\quad\, \mbox{in  }[0,T].
$$
}}
%\subsection*{Fourier's series}
% DA INSERIRE!!!

\subsection*{Well-posedness in weighted Sobolev spaces}
%Adesso introduciamo i seguenti Sobolev weighted spaces:
In order to deal with the well-posedness of problem (\ref{P2}), it is necessary to introduce the following Sobolev weighted spaces
$$\hspace{-9.0cm} H^1_a(-1,1):=$$
$$:=\{u\in L^2(-1,1): u \mbox{ locally absolutely continuous in } (-1,1), \sqrt{a} u_x \in L^2(-1,1)\}$$
and
$$ H^2_a(-1,1):=\{u\in H^1_a(-1,1)| \, au_x \in H^1 (-1,1)\}=$$
$$\hspace{-3.0cm}=\{u \in L^2 (-1,1)| u \mbox{ locally absolutely continuous in }(-1,1),$$
$$\hspace{2.5cm}au\in H^1_0 (-1,1), \, au_x \in H^1 (-1,1)\mbox{  and }
(a\,u_x)(\pm 1) %|_{x=\pm 1}
=0\} %\,(\footnote{ \mbox{
%Si può osservare che se
%One can show that if} $u\in H^2_a$ \mbox{ or } $u\in H^1_a \mbox{ then } (au_x)(t,x)|_{x=\pm 1}=0.$
%} )
$$
%Con le norme
respectively with the following norms
$$\|u\|^2_{H_a^1}:=\|u\|^2_{L^2 (-1,1)} \, + \,
|u|^2_{1,a}%\|\sqrt{a}u_x\|^2_{L^2 (-1,1)}
\mbox{ and }
\|u\|^2_{H_a^2}:=\|u\|^2_{H_a^1} \, + \, \|(au_x)_x\|^2_{L^2
(-1,1)};$$
where $|u|_{1,a}=\,\|\sqrt{a}u_x\|_{L^2 (-1,1)}$ is a seminorm.\\

%\vspace{0.5cm}
In this note we obtain the following result.
%%%%%%%%%%%%%%%%%%%%%%%%%%%%%%%%%%%%%%%%%%%%%%%%%%%%%%%%%%%%%%%%%%%%%%%%%%%%%%%%%%%%%%%

{\lemma{\label{sob2}}{ \it
\begin{equation}
H^{1}_a (-1,1)\hookrightarrow L^2(-1,1) \qquad\qquad \mbox{ with compact embedding. }
\end{equation}
}}

\proof{%(of Lemma \ref{sob2})
}{
Given $u\in H^1_a(-1,1)$, let %we consider the following extension function
$$
\bar{u}(x)=\left\{\begin{array}{l}
\displaystyle{ u \qquad\qquad\qquad\qquad \mbox{ if } x\in [-1,1]}\\ [2.5ex]
\displaystyle{ 0 \qquad\qquad\qquad\qquad \mbox{ elsewere } }~.
\end{array}\right.
$$
%By Riesz +2 Theorem
It is sufficient to prove that, for every $R>0,$
\begin{equation}\label{Riesz+2}
    \sup_{\|u\|_{1,a}\leq R}%{u\in H^1_a}
    \int_\R |\bar{u}(x+h)-\bar{u}(x)|^2\,dx\longrightarrow 0, \qquad\qquad \mbox{ as } h\rightarrow 0
\end{equation}
Let $h>0(\footnote{ In the case $h<0$ we proceed similarly.})$ and let $u\in H^1_a(-1,1)$ be such that $\|u\|_{1,a}\leq R,$ we have the following equality
$$
\hspace{-8.0cm} \int_\R |\bar{u}(x+h)-\bar{u}(x)|^2\,dx=$$
$$=\int_{-1-h}^{-1} |u(x+h)|^2\,dx\,+\int_{-1}^{1-h} |u(x+h)-u(x)|^2\,dx\,+\int_{1-h}^{1} |u(x)|^2\,dx=$$
$$=\int_{-1}^{-1+h} |u(x)|^2\,dx\,+\int_{-1}^{1-h} |u(x+h)-u(x)|^2\,dx\,+\int_{1-h}^{1} |u(x)|^2\,dx$$

First, let us prove that
%%%%%%%%%%%%%%%%%%%%%%%%%%
\begin{equation}\label{2int}
    \sup_{\|u\|_{1,a}\leq R}\int_{-1}^{1-h} |u(x+h)-u(x)|^2\,dx\longrightarrow 0, \qquad\qquad \mbox{ as } h\rightarrow 0^+.
\end{equation}
%%%%%%%%%%%%%%%%%%%%%%%%%%
Recalling that  $A(x)=\int_0^x\frac{ds}{a(s)},$ we have %let us consider
$$ |u(x+h)-u(x)|\leq \int_{x}^{x+h} \sqrt{a(s)}|u^\prime(s)|\frac{1}{\sqrt{a(s)}}\,ds \leq$$
  $$\leq \left(\int_{-1}^{1} a(s)|u^\prime(s)|^2\,ds\right)^{\frac{1}{2}}\left(\int_{x}^{x+h} \frac{ds}{a(s)}\,\right)^{\frac{1}{2}}= |u|_{1,a}\left[A(x+h)-A(x)\right]^{\frac{1}{2}}.$$
By integrating on $[-1,1-h]$, since $A\in L^1 (-1,1)$ (by assumption 3.b)),  we obtain
$$\int_{-1}^{1-h} |u(x+h)-u(x)|^2\,dx\leq |u|^2_{1,a}\int_{-1}^{1-h}\left(A(x+h)-A(x)\right)\, dx \leq$$
$$\hspace{-5cm}\leq %|u|^2_{1,a}
R^2\left[\int_{-1+h}^{1} A(x)\, dx\,-\,\int_{-1}^{1-h} A(x)\, dx\,\right]=$$
$$\hspace{2.5cm}=%|u|^2_{1,a}
R^2\left[\int_{1-h}^{1} A(x)\, dx\,-\,\int_{-1}^{-1+h} A(x)\, dx\,\right]\longrightarrow 0,\quad \mbox{ as } h\rightarrow 0^+.$$

%%%%%%%%%%%%%%%%%%%%%%%
Now, let us prove that
\begin{equation}\label{3int}
     \sup_{\|u\|_{1,a}\leq R}\int_{1-h}^{1} |u(x)|^2\,dx\longrightarrow 0, \qquad\qquad \mbox{ as } h\rightarrow 0^+.
\end{equation}
%%%%%%%%%%%%%%%%%%%%%%%%%%%%%%%%%%%%%%%%%%%%%%%%%%%%%%%%%%%%%%%
%%%%%%%%%Teor. fondam. calcolo integrale%%%%%%%%%%%%%%%%%%%%%%%%%%%%%%%%%%%%%%%
%\noindent We start to consider the following equality
%\begin{equation}\label{ffond}
%u(x)-u(0)=\int_{0}^x u^\prime(s)\,ds.
%\end{equation}
%%%%%%%%%%%%%%%%%%%%%%%%%%%%%%%%%%%%%%%%%%%%%%%%%%%%%%%%%%%%%%%%%%%%%%%%%%%%%%%%
We have
$$|u(0)|\leq |u(x)|+\int_{0}^{x} \sqrt{a(s)}|u^\prime(s)|\frac{1}{\sqrt{a(s)}}\,ds\leq$$
$$\leq |u(x)|+\left(\int_{-1}^{1} a(s)|u^\prime(s)|^2\,ds\right)^{\frac{1}{2}}\left(\int_{0}^{x} \frac{ds}{a(s)}\,\right)^{\frac{1}{2}}\leq |u(x)|+ |u|_{1,a}\sqrt{A(x)}\,.$$

\noindent By integrating on $[0,1],$ we obtain
$$\hspace{-2cm}|u(0)|\leq \int_{0}^{1}|u(x)|\, dx+|u|_{1,a}\int_{0}^{1}\sqrt{A(x)}\,dx\leq$$
$$\hspace{2cm}\leq\|u\|_{L^2(-1,1)}+|u|_{1,a}\int_{0}^{1}\sqrt{A(x)}\,dx\leq C\|u\|_{1,a}\,.$$
Then,
\begin{equation}\label{u0}
  |u(0)|\leq C\,R%\|u\|_{1,a}
  \,.
\end{equation}
Now, %by (\ref{ffond}) and (\ref{u0}),
it follows that
$$|u(x)|^2\leq 2|u(0)|^2+2 A(x)|u|^2_{1,a}\leq C\,R^2%\|u\|^2_{1,a}
+2 A(x) R^2%|u|^2_{1,a}
\,.$$
Finally, since %because
$A\in L^1(-1,1)$, by integrating on $[1-h,1]$ we obtain
$$\int_{1-h}^{1}|u(x)|^2\, dx\leq C\,h R^2%\|u\|^2_{1,a}
+2 %|u|^2_{1,a}
R^2\,\int_{1-h}^{1}A(x)\,dx \longrightarrow 0, \qquad\qquad \mbox{ as } h\rightarrow 0^+. $$
%%%%%%%%%%%%%%%%%%%%%%%

Similarly, we can prove that
\begin{equation}\label{1int}
    \sup_{\|u\|_{1,a}\leq R}\int_{-1}^{-1+h} |u(x)|^2\,dx\longrightarrow 0, \qquad\qquad \mbox{ as } h\rightarrow 0^+.
\end{equation}

%%%%%%%%%%%%%%%%%%%%%%%%%%%%%%%%%%%%%%%%%%%%%%%%%%%%%%%%%%%%%%%
\noindent By (\ref{2int}), (\ref{3int}) and (\ref{1int}) we obtain (\ref{Riesz+2}).\hspace{4cm} \eproof}\\
%%%%%%%%%%%%%%%%%%%%%%%%%%%%%%%%%%%%%%%%%%%%%%%%%%%%%%%%%%%%%%%
%%%%%%%%%%%%%%%%%%%%%%%%%%%%%%%%%%%%%%%%%%%%%%%%%%%%%%%%%%%%%%%%%%%%%%%%%%%%%%%%%%%%%%%

\vspace{0.5cm}
We now recall the existence and uniqueness result for system (\ref{P2}) obtained in \cite{CMP} (see also \cite{ACF}).
 %The well-posedness of system (\ref{P2}),
%%è assicurata da un risultato contenuto
%is shown in \cite{ACF} and in \cite{CMP}.
%%\begin{equation} \label{1.3} D(A) \, = \, H^2_a \,\, \mbox{and}
%%\,\, \forall \, u \, \in \, D(A), \,\,\, Au:=(au_x)_x.
%%\end{equation}
%%Infatti, nelle referenze sopra citate si definisce l'operatore
%In fact, in the above mentioned references,
Let us consider, first, the operator $(A_0,D(A_0))$
%mediante
defined by
 \begin{equation}\label{D(A0)}
   \left\{\begin{array}{l}
\displaystyle{D(A_0)=H^2_a (-1,1) }\\ [2.5ex]
\displaystyle{A_0u=(au_x)_x\,, \,%\qquad\qquad\qquad\qquad\quad
\,\, \forall \,u \in D(A_0)}\,.
\end{array}\right.
 \end{equation}
Observe that $A_0$ is a closed, self-adjoint, dissipative operator with dense domain in $L^2 (-1,1)$.
Therefore, $A_0$ is the infinitesimal generator of a $C_0-\mbox{semigroup}$ of contractions in $L^2 (-1,1)$.

Next, given $\alpha\in L^\infty (-1,1),$ let us introduce the operator

 \begin{equation}\label{D(A)}
 \left\{\begin{array}{l}
\displaystyle{D(A)=D(A_0) }\\ [2.5ex]
\displaystyle{A = A_0 + \alpha I\, %\qquad\qquad\qquad\qquad\quad
%\,\, \forall \,u \in D(A_0)
}\,.
\end{array}\right.
 \end{equation}
 %Then, we have the following
  %where the operator $A_0$ is defined in (\ref{D(A0)}) and $\alpha\in L^\infty (-1,1),$
%is a bounded perturbation of $A_0.$
%  $D(A)=H^2_a$ and $ Au:=(au_x)_x,\,\forall \,u \in D(A) .$\\
%per tale operatore si ha la seguente proposizione:
For such an operator we have the following proposition.

%{\proposition{\label{A_0}}{\it
%The operator $A_0: D(A_0)\longrightarrow L^2(-1,1)$ is closed,
%self-adjoint and negative with dense domain, $H^{2}_a (-1,1)\hookrightarrow L^2(-1,1)\, \mbox{ with compact embedding. %}$
%}}
%\vspace{0.5cm}
{\proposition{\label{A_0}}{\it
\begin{itemize}
\item $D(A)$ is compactly embedded and dense in $L^2(-1,1)$.
\item $A: D(A)\longrightarrow L^2(-1,1)$ is %a closed,
%self-adjoint and strictly dissipative operator.
the infinitesimal generator of a strongly continuous semigroup, $e^{tA}$, of bounded linear operators on $L^2(-1,1)$.
\end{itemize}
}}
\vspace{0.3cm}
Observe that problem (\ref{P2}) can be recast in the Hilbert space $L^2(-1,1)$ as
\begin{equation}\label{Ball}
 \left\{\begin{array}{l}
\displaystyle{u^\prime(t)=A\,u(t)\,,\qquad  t>0 }\\ [2.5ex]
\displaystyle{u(0)=u_0\, %\qquad\qquad\qquad\qquad\quad
%\,\, \forall \,u \in D(A_0)
}~.
\end{array}\right.
\end{equation}
where $A$ is the operator in (\ref{D(A)}).

\vspace{0.5cm}
We recall that a \textit{weak solution} of (\ref{Ball}) is a function $u\in C^0([0,T];L^2(-1,1))$ such that for every $v\in D(A^*)$ the function $\langle u(t),v\rangle$ is absolutely continuous on $[0,T]$
and
$$\frac{d}{dt}\langle u(t),v\rangle=\langle u(t),A^*v\rangle\,,$$
for almost $t\in [0,T]$
 (see \cite{Ba}).
%Hence, $A$ is the infinitesimal generator of a strongly continuous
%semigroup $e^{tA}$ on $L^2(-1,1)$.
%%Moreover, since $A_0$ is a
%%generator and the operator $B(t)$ defined as
%%$$
%%B(t)u:=-c(t,\cdot)u
%%$$
%%can be seen as a bounded perturbation of $A_0$.

\vspace{0.5cm}
%Working in the spaces considered above, we have that (\ref{P2}) is well-posed in the sense of semigroup theory.

{\theorem{}{\it
For every $\alpha\in L^\infty((0,T)\times(-1,1))$
and every $u_0 \in L^2(-1,1)$, there exists a unique
$$u\in C^0([0,T];L^2(-1,1))\cap L^2(0,T;H^1_a (-1,1))$$
weak solution to (\ref{P2}), which coincides with $e^{tA}u_0.$
}}

\vspace{1cm}
%Definiamo la seguente norma nello spazio
In the space
$$\mathcal{B}(0,T)=C^0([0,T];L^2(-1,1))\cap L^2(0,T;H^1_a (-1,1))$$
let us define the following norm
\begin{equation}
\label{normaB}
 \|u\|^2_{\mathcal{B}(0,T)}= \sup_{t\in
[0,T]}\|u(t,\cdot)\|^2_{L^2(-1,1)}+2\int^T_{0}\int^1_{-1}a(x)u^2_x
dx\,,\,\,\forall u \in \mathcal{B}(0,T)\,.
\end{equation}

\section{Some auxiliary lemmas and the proofs of main results}

Let $A = A_0 + \alpha I,$ where the operator $A_0$ is defined in (\ref{D(A0)}) and $\alpha\in L^\infty (-1,1).$
Since $A$ is self-adjoint and $D(A)\hookrightarrow L^2(-1,1)$ is compact (see Proposition \ref{A_0}), we have the following (see also \cite{BR}).

%\vspace{-0.5cm}

{\lemma{\label{spectrum}}{\it
 %is closed,
%self-adjoint and negative with dense domain, $H^{2}_a (-1,1)\hookrightarrow L^2(-1,1)\, \mbox{ with compact embedding.}$
%Then
There exists an increasing sequence $\{\lambda_k\}_{k\in\N},$ with 
$\lambda_k\longrightarrow +\infty\, \mbox{ as } \, k \, \rightarrow\infty\,,$
such that the eigenvalues of $A$ are given by $\{-\lambda_k\}_{k\in\N}$, and the corresponding eigenfunctions $\{\omega_k\}_{k\in\N}$ form a complete orthonormal system in $L^2(-1,1)$.
%$$\{-\lambda_k\}_{k\in\N}\qquad\qquad\qquad\mbox{  and  }\qquad\qquad\qquad\{\omega_k\}_{k\in\N},$$
%respectively the eigenvalues and
%the orthonormalized in $L^2(-1,1)$ eigenfunctions of the operator $A,$
%%spectral problem
%%$$A\omega=\lambda\omega,\qquad\qquad \omega\in H^{2}_a (-1,1),$$
%where $\{\lambda_k\}$ is a increasing sequence such that
%$$\lambda_k\longrightarrow +\infty\qquad \mbox{ as } \, k \, \rightarrow\infty\,.$$
}}

 %%Quindi in definitiva si può considerare
%So we can consider
\vspace{0.5cm}
In this note we obtain the following result

{\lemma{\label{Autof}}}{\it
Let $v\in C^\infty([-1,1]), v >0$ on $[-1,1],$ let $\alpha_*(x)=-\frac{(a(x)v_{x}(x))_x}{v(x)},\, x\in
(-1,1).$ Let A be the operator defined in (\ref{D(A)}) with $\alpha=\alpha_*$

\begin{equation}\label{operalfastella}
%\left\{\begin{array}{l}
%\displaystyle{v_t-(a(x) v_x)_x =\alpha_* (x)v\,\,\qquad \mbox{in} \qquad Q_T \,=\,(0,T)\times(-1,1) }\\ [2.5ex]
%\displaystyle{a(x)v_x(t,x)|_{x=\pm 1} = 0\,\,\qquad\qquad\qquad\,\, t\in (0,T) \qquad\qquad\qquad\qquad %(\ref{P2})
% }\\ [2.5ex]
%%\displaystyle{v(0,x)=v_0 (x) \,\qquad\qquad\qquad\qquad\,\,\,\,\, x\in (-1,1)}~,
%\end{array}\right.
\left\{\begin{array}{l}
\displaystyle{D(A)=H^2_a (-1,1)}\\ [2.5ex]
\displaystyle{A = A_0 + \alpha_* I}~,
%\displaystyle{v(0,x)=v_0 (x) \,\qquad\qquad\qquad\qquad\,\,\,\,\, x\in (-1,1)}~,
\end{array}\right.
\end{equation}
%where the operator $A_0$ is defined in (\ref{D(A0)}).
and let $\{\lambda_k\}, \{\omega_k\}$ be the eigenvalues and eigenfunctions of $A,$ respectively, given by Lemma \ref{spectrum}.
Then $$\lambda_1=0\,\, \mbox{ and } \,\,|\omega_1|=\frac{v}{\|v\|%_{L^2(-1,1)}
}.$$
Moreover, $\frac{v}{\|v\|%_{L^2(-1,1)}
}$ and $-\frac{v}{\|v\|%_{L^2(-1,1)}
}$ are the only eigenfunctions of $A$ with norm $1$ that do not change sign in $(-1,1)$.

%The function
%$$\frac{v}{\|v\|%_{L^2(-1,1)}
%}$$
%is the first eigenfunction of the (\ref{operalfastella}) operator and it is associated to zero-eigenvalue. Moreover, %this is the only orthonormal eigenfunction of (\ref{operalfastella})
%%che non cambia segno in
%that doesn't change sign in $%\Omega=
%(-1,1)$.
}

%\left\{\begin{array}{l}
%\displaystyle{v_t-(a(x) v_x)_x =\alpha_* (x)v\,\,\qquad \mbox{in} \qquad Q_T \,=\,(0,T)\times(-1,1) }\\ [2.5ex]
%\displaystyle{a(x)v_x(t,x)|_{x=\pm 1} = 0\,\,\qquad\qquad\qquad\,\, t\in (0,T) \qquad\qquad\qquad\qquad %(\ref{P2})
% }\\ [2.5ex]
%%\displaystyle{v(0,x)=v_0 (x) \,\qquad\qquad\qquad\qquad\,\,\,\,\, x\in (-1,1)}~,
%\end{array}\right.

{\remark{Problem (\ref{operalfastella}) is equivalent to the following differential problem
\begin{equation}\label{stella}
\left\{\begin{array}{l}
\displaystyle{(a(x) \omega_x)_x +\alpha_* (x)\omega+\lambda\,\omega=0\,\,\qquad \mbox{in} \qquad (-1,1) %\,=\,(0,T)\times(-1,1)
}\\ [2.5ex]
\displaystyle{a(x)\omega_x(x)|_{x=\pm 1} = 0\,\,%\qquad\qquad\qquad\,\, %t\in (0,T) \qquad\qquad\qquad\qquad %(\ref{P2})
\,\,\, }~.
%\displaystyle{v(0,x)=v_0 (x) \,\qquad\qquad\qquad\qquad\,\,\,\,\, x\in (-1,1)}~,
\end{array}\right.
\end{equation}}}

\proof{ (of Lemma \ref{Autof})}{
\underline{STEP.1}
%Scelto un qualsiasi stato iniziale nonnegative, non nullo
%e un qualsiasi
%and any target state $v_d$
%%come descritto in
%as described in (\ref{vd}) in STEP.1,
%%poniamo
%we set
%\begin{equation}\label{2.9}
%\alpha_*(x)=-\frac{(a(x)v_{dx}(x))_x}{v_d(x)},\qquad \, x\in
%(-1,1).
%\end{equation}
%%Si può osservare che
%Then, by (\ref{vd}),
%$$\alpha_* (x)\in L^\infty(-1,1)\,.  %  \mbox{ and }  \alpha_* (x)\not\equiv0
%$$
%%in fact if $\alpha_*(x)\equiv 0, \mbox{
%%then } a(x)v_{dx}(x)=$
%%costante e poichè
%%constant and since $\lim_{x\longrightarrow\pm 1}a(x)v_{dx}(x)=0$
%%allora......
%%then --
%%LO DEVO SISTEMARE!!! %GFAgosto2010
%%......e quindi
%%-- and therefore $v_d(x)=0$,
%%che è in contrasto con la
%%that is in contrast with (\ref{vd}).
%%\\
%%%non è uguale a zero
%%%in $L^\infty (-1,1) (\mbox{diversamente se }
%%%y_{dxx}(x\equiv0)\mbox{ e quindi essendo } y_d(-1)=y_d(1)=0,
%%\mbox{ ne seguirebbe}\, y_d\equiv 0\quad \mbox{in  }\Omega)$. Gli
%%%autovalori associati con $\alpha_*$ ci sono indicati inoltre da
%%%$\{\lambda_k\}^\infty_{k=1}$.
%%%\\ La \ref{2.9} vuole dire che la funzione
%%%è un'autofunzione per
%%Gli autovalori associati al problema
We denote by
$$\{-\lambda_k\}_{k\in\N}\qquad\qquad\qquad\mbox{  and  }\qquad\qquad\qquad\{\omega_k\}_{k\in\N},$$
respectively, the eigenvalues and
%the
orthonormal eigenfunctions of the
%spectral problem $A\omega=\lambda \omega,$ with $A=A_0+\alpha_*I$
operator (\ref{operalfastella}) (see Lemma \ref{spectrum}).
Therefore,
\begin{equation}\label{3}
\langle\omega_k,\omega_h\rangle_{L^2(-1,1)}=\int^1_{-1}
\omega_k(x)\omega_h(x)dx=0, \qquad \mbox{ if }h\neq k\,.
\end{equation}
We can see, by easy calculations, that an eigenfunction of the operator defined in (\ref{operalfastella})
%è la funzione:
is the function
$$
\frac{v(x)}{\|v\|%_{L^2(\Omega)}
}\,,
$$
%%%%%%%%%%%%%%%%%%%%%%%%%%%%%%%%%%%%%%%%%%%%%%%%%%%%%%%%%%%%%%%%%%%%%%%%%%%%%%%%%%%%%%
%%  VERIFICA CHE \alpha_* è un'autofunzione
%%%%%%%%%%%%%%%%%%%%%%%%%%%%%%%%%%%%%%%%%%%%%%%%%%%%%%%%%%%%%%%%%%%%%%%%%%%%%%%%%%%%%%
%Infatti
%In fact, we have
%$$
%\bigg(a(x)\frac{v_{dx}(x)}{\|v_d\|_{L^2(\Omega)}}\bigg)_x +
%\alpha_* \frac{v_d(x)}{\|v_d\|_{L^2(\Omega)}} - \lambda
%\frac{v_d(x)}{\|v_d\|_{L^2(\Omega)}} =
%$$
%$$
%=\frac{(a(x)v_{dx}(x))_x}{\|v_d\|_{L^2(\Omega)}}-\frac{(a(x)v_{dx}(x))_x}{v_d(x)}\frac{v_d(x)}{\|v_d\|_{L^2(\Omega)}}-\lambda
%\frac{v_d(x)}{\|v_d\|_{L^2(\Omega)}}=0\,,
%$$
%da cui
%from which it follows that $\frac{v_d(x)}{\|v_d\|_{L^2(\Omega)}}$
%%è un'autofunzione per
%is an eigenfunction of (\ref{1.3})
%associata all'autovalore
associated with the eigenvalue $\lambda=0$.
%%Indichiamo con
%We denote by $\{\omega_k\}^\infty_{k=1}$
%%le autofunzioni ortogonali in
%the orthonormal eigenfunctions in $L^2(\Omega)$
%%del problema
%of the problem (\ref{1.3}),
%%quindi
%Per quanto detto e considerato in particolare che
%%%%%%%%%%%%%%%%%%%%%%%%%%%%%%%%%%%%%%%%%%%%%%%%%%%%%%%%%%%%%%%%%%%%%%%%%%%%%%
%% FINE VERIFICA AUTOFUNZIONE
%%%%%%%%%%%%%%%%%%%%%%%%%%%%%%%%%%%%%%%%%%%%%%%%%%%%%%%%%%%%%%%%%%%%%%%%%%%%%%
Taking into account the above and considering that
$v(x)>0,\,\forall x \in (-1,1)$
\begin{equation}
\label{4}
\,\exists\, k_* \in \N\,\,:\omega_{k_*}(x)=\frac{v(x)}{\|v\|%_{L^2(\Omega)}
}>0\,  \mbox{ or }  \,\omega_{k_*}(x)=-\frac{v(x)}{\|v\|%_{L^2(\Omega)}
}<0,\,\,\forall x\in (-1,1)\,.
\end{equation}
%Scritta la
Writing (\ref{3}) with $k=k_*$
%si ha:
 we obtain
\begin{equation}
\label{5} \langle\omega_{k_*},\omega_h\rangle_{L^2(-1,1)}=\int^1_{-1}
\omega_{k_*}(x)\omega_h(x)dx=0, \qquad \forall h\neq k_*\, .
\end{equation}
%Quindi considerando la
Therefore, considering (\ref{5})
%e tenuto conto che
and keeping in mind that $\omega_{k_*}>0$ or $\omega_{k_*}<0$ in $(-1,1)$,
%si ha che
we observe that $\omega_{k_*}$
%è l'unica autofunzione ortonormale di
is the only eigenfunction of the operator defined in (\ref{operalfastella})
%che non cambia segno in
that doesn't change sign in $%\Omega=
(-1,1)$.\\
%Inoltre dal fatto che
%Since $v_0\in L^2(-1,1)$ %, \,v_0\geq 0 \mbox{ and } v_0\not\equiv0$ in $(-1,1),$
%%si ha
%and $\int^1_{-1} v_0(x)v_d(x)dx>0$
%we obtain
%\begin{equation}
%\label{6}
%\int^1_{-1} v_0(x)\omega_{k_*}(x)dx>0.
%\end{equation}

\noindent
\underline{STEP.2}
%Proviamo adesso che
Let us now prove that
\begin{equation}
\label{6bis} k_*=1\,,
%\lambda_1=0
\end{equation}
%%e quindi di conseguenza
%and consequently
that is, $\lambda_1=0$.
By a well-known variational characterization of the first eigenvalue, we have %the following representation
$$\lambda_1=\inf_{u\in H^1_a (-1,1)}\frac{\int^1_{-1} \left(a\,u_x^2\,-\alpha_*\,u^2\right)\, dx}{\int^1_{-1} u^2\, dx}\,\,.$$
By Lemma \ref{spectrum}, since $\lambda_{k_*}= 0,$ it is sufficient to prove that $\lambda_1\geq 0$, or\begin{equation}\label{VARIAT}
  \int^1_{-1}\alpha_*\,u^2\, dx\leq   \int^1_{-1} a\,u_x^2\,dx,\qquad \forall\,u\in H^1_a(-1,1)
\end{equation}
Integrating by parts, we have
$$
\int^1_{-1}\alpha_*\,u^2\, dx=-\int^1_{-1}\frac{(a\,v_x)_x}{v}\,u^2\, dx = \int^1_{-1}a\,v_x\left(\frac{u^2}{v}\right)_x\,dx=
$$
$$
=\int^1_{-1}a\,v_x\frac{2 u u_x}{v}\,dx-\int^1_{-1}a\,v^2_x\left(\frac{u^2}{v^2}\right)\,dx=
$$
$$
=2\int^1_{-1}\sqrt{a}\,\frac{v_x}{v}u\sqrt{a}u_x\,\, dx-\int^1_{-1}a\,v^2_x\left(\frac{u^2}{v^2}\right)\,\, dx\leq
$$
$$\leq\int^1_{-1}\,a\,\left(\frac{ v_x u}{v}\right)^2\,dx+\int^1_{-1} a u^2_x\,dx-\int^1_{-1}a\,v^2_x\left(\frac{u^2}{v^2}\right)\,dx = \int^1_{-1} a u^2_x\,\, dx\,,$$
from which (\ref{VARIAT}).
\hfill
\eproof}

\vspace{0.5cm}
For the proof of  Theorem \ref{T1}
%è necessario il seguente
the following Lemma %\ref{NN}
 is necessary.

% \vspace{-0.5cm}
%\newpage
{\lemma{\label{NN}}{\it
%Sia
Let $T>0$,
%, sia
$\alpha\in L^\infty(Q_T)$, let $v_0 \in L^2(-1,1),\,v_0(x)\geq 0 \,\mbox{ a.e. } x \in (-1,1)$ and % $\alpha(x)\geq 0,$
%sia
let $v\in\mathcal{B}(0,T)$ %= C([0,T],\, L^2(-1,1))\, \cap \,L^2([0,T],H^1_a (-1,1))$
%soluzione del sistema lineare
 be the solution to the linear system
%\begin{equation*}%\label{P2bis}
$$\left\{\begin{array}{l}
\displaystyle{v_t-(a(x) v_x)_x =\alpha (t,x)v\,\,\qquad \mbox{in} \qquad Q_T \,=\,(0,T)\times(-1,1) }\\ [2.5ex]
\displaystyle{a(x)v_x(t,x)|_{x=\pm 1} = 0\,\,\qquad\qquad\qquad\,\, t\in (0,T) \qquad\qquad\qquad\qquad %(\ref{P2})
 }\\ [2.5ex]
\displaystyle{v(0,x)=v_0 (x) \,\qquad\qquad\qquad\qquad\,\,\,\,\, x\in (-1,1)}~.
\end{array}\right.$$
%\left\{\begin{array}{l}
%{ \displaystyle{ v_t-(a(x) v_x)_x - \alpha v=0 }}\\[2.5ex]
%{ \displaystyle{ a(x)v_x(t,x)|_{x=\pm 1} = 0 }}\\[2.5ex]
%{ \displaystyle{ v(0,x)=v_0 (x) }}~.
%\end{array}\right.
%%%
%\begin{cases}
%v_t-(a(x) v_x)_x - \alpha v=0
%\\
%a(x)v_x(t,x)|_{x=\pm 1} = 0
%\\
%v(0,x)=v_0 (x)
%\end{cases}
%\end{equation*}
Then\\
$$v(t,x)\geq 0,\,\,\,\,\forall (t,x)\in Q_T\,.%>0
$$
%%%%%%%%%%%%%%%
%Moreover, if $v_0 \in L^\infty (-1,1),\,v_0(x)\geq 0 \,\mbox{ a.e. } x \in (-1,1),$ we have also the following %estimate
%$$0\leq v(t,x)\leq e^{\|\alpha\|_\infty
%t}\|v_0\|_{L^\infty (-1,1)},\,\,\,\,t\in [0,T]\,.%>0
%$$
}}

%%%%%%%%%%%%%%%

\proof{
%%\begin{itemize}
%%\item
%\underline{STEP 1} \textit{
%%Proviamo che
%Let us prove that $v(t,x)\geq 0 \,\, \forall(t,x)\in Q_T$.}\\
%%Sia
Let $v\in\mathcal{B}(0,T)$ %$v\in $$\cal{B}$$(0,T)= C([0,T],\, L^2(-1,1))\, \cap \,L^2([0,T],H^1_a (-1,1))$
be the
%soluzione del sistema
solution to the system (\ref{P2}),
%consideramo le seguenti funzioni:
and we consider the positive-part
%$$v^+\,=\,\max (v,0)\,=
%\left\{\begin{array}{l}
%\displaystyle{ v\qquad\qquad \mbox{ in }\{v>0\}} \\ [2.5ex]
%\displaystyle{ 0\qquad\qquad \mbox{ in }\{v\leq 0 \}}~.
%\end{array}\right.$$
and the negative-part.
%%\,\begin{cases}
%% v\qquad\qquad \mbox{in}\{v>0\} \\
%%0\qquad\qquad \mbox{in}\{v\leq 0\}
%%\end{cases}
%$$v^-\,=\,\max (0,-v)\,=\,
%\left\{\begin{array}{l}
%\displaystyle{-v\qquad\,\,\,\, \mbox{ in }\{v<0\}}\\ [2.5ex]
%\displaystyle{ 0\qquad\qquad \mbox{ in }\{v\geq 0\} }~.
%\end{array}\right.$$
%%\begin{cases}
%% -v\qquad\,\,\,\, \mbox{in}\{v<0\}\\
%%0\qquad\qquad \mbox{in  }\{v\geq 0\}
%%\end{cases}$$
%%si ha, come richiamato nei preliminari
%%Therefore we have, as recalled in the preliminaries (Section 2)
%%$
%Since $v=v^+-v^-\,,$
%Osservato che
%Having observed that
%Since $v^+,\,v^-\geq 0$,
%per provare la tesi è sufficiente provare che
 %to prove the thesis,
It is sufficient to prove that
$$v^-(t,x)\equiv 0\qquad\quad\mbox{in  }Q_T\,.$$
%Moltiplicando ambo i membri dell'equazione del problema per
Multiplying both members equation of the problem (\ref{P2}) by
$v^-$ %ed integrando in
and integrating it on $(-1,1)$
%si ha:
we obtain
\begin{equation}
\label{4.2}
\int^1_{-1}\left[ v_t v^- -(a(x)v_x)_x v^- -\alpha v
v^- \right]dx=0.
\end{equation}
%Tenendo conto della definizione di si ha:
Recalling the definition $v^+$ and $v^-$, we
obtain
$$
\int^1_{-1}v_t v^- dx = \int^1_{-1} (v^+ - v^-)_t v^- dx =
-\int^1_{-1} (v^-)_t v^- dx = -\frac{1}{2} \frac{d}{dt} \int
(v^-)^2 dx\,.
$$
%Inoltre, integrando per parti, e applicando il si ha
Integrating by parts and applying  Theorem \ref{A.1}, we
obtain $v^- \in H^1_a (-1,1)$
%e vale la seguente uguaglianza
and the following equality
$$
\int^1_{-1} (a(x)v_x)_x v^-\,dx =[a(x)v_x v^-]^1_{-1} - \int^1_{-1}
a(x)v_x(-v)_x\,dx = \int^1_{-1} a(x) v^2_x\,dx\,.
$$
%Si ha anche
We also have
$$
\int^1_{-1}\alpha v v^- dx = -\int^1_{-1}\alpha(v^-)^2 dx
$$
%Quindi la  diventa:
and therefore (\ref{4.2}) becomes
$$ -\frac{1}{2} \frac{d}{dt}
\int^1_{-1}(v^-)^2 dx + \int^1_{-1}\alpha (v^-)^2 dx = \int^1_{-1}
a(x) v^2_x\geq 0,
$$
%da cui,
from which
$$
\frac{d}{dt}\int^1_{-1}(v^-)^2 dx\leq 2 \int^1_{-1}\alpha (v^-)^2
dx\leq 2\|\alpha\|_\infty\int^1_{-1}(v^-)^2 dx.
$$
%da ciò applicando il
From the above inequality, applying Gronwall's lemma %( Lemma \ref{GW} with $\eta(t)=\int^1_{-1}(v^-)^2 dx,\,
%\phi(t)=2\|\alpha\|_\infty,\, \psi(t)\equiv 0$)
%si ha:
 we obtain
$$
\int^1_{-1}(v^- (t,x))^2 dx\leq
e^{2t\|\alpha\|_\infty}\,\,\int^1_{-1}(v^-(0,x))^2dx\,.%=0\,.
$$
%in quanto essendo
Since
$$v(0,x)=v_0(x)\geq 0\,,$$
%segue
we have
$$v^-(0,x)=0.$$
%Quindi abbiamo provato che
Therefore,% we prove that
$$
v^-(t,x)=0, \qquad\qquad \,\,\,\quad\forall (t,x)\in Q_T.
$$
%Da ciò segue che
From this, as we mentioned initially, it follows that
$$
v(t,x)=v^+(t,x)\geq 0 \qquad\forall (t,x)\in Q_T.
$$
\hspace{12cm}\eproof}
We are now ready to prove our main result.\\

\proof{ (of Theorem \ref{T1})}{

\noindent
\underline{STEP.1}
 %Per provare il
 To prove Theorem \ref{T1}
 %è sufficiente considerare un qualsiasi insieme di non-negativi
 it is sufficient to consider the set of target states %$v_d$
 %%denso nell'isieme di tutti i non negativi elementi di
 %dense in the set of all the nonnegative elements of $L^2(-1,1)$.\\
 %%Per la densità di
 %Because of the density of $C^{\infty}(-1,1)$ in $L^2(-1,1)$,
 %%per ogni non negativa funzione
 %for every nonnegative function
\begin{equation}\label{vd}
 \qquad v_d\in C^\infty ([-1,1]), %H^2(-1,1)\cap %H^1_0(-1,1) \mbox{  such that  } v_d (x)>0 \mbox{  in  } (-1,1) \mbox{  and  } \frac{(a(x)v_d(x))_x}{v_d(x)}\in L^\infty(-1,1),
  \qquad\,v_d>0  \mbox{  on  } [-1,1].
 \end{equation}
 Indeed, regularizing by convolution, every function $v_d\in L^2(-1,1), v_d\geq 0$ can be approximated by a sequence of strictly positive $C^\infty ([-1,1])-$ functions.\\

 %%%%%%%%%%%%%%%%%%%%%%%%%%%%%%%%%%%%%%%%%%%%%%%%%%%%%%%%%%%%%%%%%%%%%%%%%%%%%%%%%%
 %APPROSSIMAZIONE CON DETTAGLI
 %%%%%%%%%%%%%%%%%%%%%%%%%%%%%%%%%%%%%%%%%%%%%%%%%%%%%%%%%%%%%%%%%%%%%%%%%%%%%%%%%%
% In fact, at first every function $v_d\in L^2(-1,1), v_d\geq 0$ can be approximated by the following sequence of %strictly positive $L^2(-1,1)-$functions
%$$\{\tilde{v}_{dk}\}_{k\in\N}=\left\{v_{d}+\frac{1}{k}\right\}_{k\in\N}\subseteq L^2(-1,1) $$
 %%tale che
 %so that
 %$$v_{dk}\stackrel{L^2}{\longrightarrow} v_d, \mbox{ as } k\rightarrow\infty .$$
 %%Inoltre ciascua funzione
 %After, regularizing by convolution, $L^2(\R)$-extension function of every $\tilde{v}_{dk}\in L^2(-1,1)$ can be %approximated in $L^2(\R)$ by a sequence of strictly positive $C^\infty(\R)-$functions %$\{\tilde{v}^\varepsilon_{dk}\}_{\varepsilon>0}$
 %so that
 %$$\tilde{v}^\varepsilon_{dk}\stackrel{L^2}{\longrightarrow} v_d, \mbox{ as } k\rightarrow\infty \mbox{ and } %\varepsilon\rightarrow 0^+ .$$
 %%%si può approssimare (andando a smussare gli eventuali punti angolosi) con funzioni di classe
 %%can be approximate (smoothesing any angulr points) with functions from class $C^{\infty}(-1,1)$
 %%%con al più un numero finito d punti di discontinuità di prima specie per la derivata seconda.
%%%with a finite number on discontinuous points of the first kind for
%%%the deriving second.
 %%%%%%%%%%%%%%%%%%%%%%%%%%%%%%%%%%%%%%%%%%%%%%%%%%%%%%%%%%%%%%%%%
 %Fine approssimazione con dettagli
 %%%%%%%%%%%%%%%%%%%%%%%%%%%%%%%%%%%%%%%%%%%%%%%%%%%%%%%%%%%%%%%%%

 %%Quindi in definitiva si può considerare
%So we can consider

\smallskip\noindent
\underline{STEP.2}
%Scelto un qualsiasi stato iniziale nonnegative, non nullo
Taking any nonzero, nonnegative initial state $v_0\in
L^2(-1,1)$
%e un qualsiasi
and any target state $v_d$
%come descritto in
as described in (\ref{vd}) in STEP.1,
%poniamo
let us set
\begin{equation}\label{2.9}
\alpha_*(x)=-\frac{(a(x)v_{dx}(x))_x}{v_d(x)},\qquad \, x\in
(-1,1).
\end{equation}
%Si può osservare che
Then, by (\ref{vd}),
$$\alpha_* (x)\in L^\infty(-1,1)\,.  %  \mbox{ and }  \alpha_* (x)\not\equiv0
$$
%in fact if $\alpha_*(x)\equiv 0, \mbox{
%then } a(x)v_{dx}(x)=$
%costante e poichè
%constant and since $\lim_{x\longrightarrow\pm 1}a(x)v_{dx}(x)=0$
%allora......
%then --
%LO DEVO SISTEMARE!!! %GFAgosto2010
%......e quindi
%-- and therefore $v_d(x)=0$,
%che è in contrasto con la
%that is in contrast with (\ref{vd}).
%\\
%%non è uguale a zero
%%in $L^\infty (-1,1) (\mbox{diversamente se }
%%y_{dxx}(x\equiv0)\mbox{ e quindi essendo } y_d(-1)=y_d(1)=0,
%\mbox{ ne seguirebbe}\, y_d\equiv 0\quad \mbox{in  }\Omega)$. Gli
%%autovalori associati con $\alpha_*$ ci sono indicati inoltre da
%%$\{\lambda_k\}^\infty_{k=1}$.
We denote by
$$\{-\lambda_k\}_{k\in\N}\qquad\qquad\qquad\mbox{  and  }\qquad\qquad\qquad\{\omega_k\}_{k\in\N},$$
respectively, the eigenvalues and orthonormal eigenfunctions\footnote{As first eigenfunction we take the one which is positive in $(-1,1)$.} of the
spectral problem $A\omega+\lambda \omega=0,$ with $A=A_0+\alpha_*I$ (see Lemma \ref{spectrum}).\\
%Therefore
%\begin{equation}
%\label{3} \langle\omega_k,\omega_h\rangle_{L^2(-1,1)}=\int^1_{-1}
%\omega_k(x)\omega_h(x)dx=0, \qquad \mbox{ if }h\neq k\,.
%\end{equation}

%%The eigenvalues associated to the operator %(\ref{1.3})
%%%li indichiamo con
%%are denoted by $\{-\lambda_k\}^\infty_{k=1}.$
%%%Osserviamo che un'autofunzione dell'equazione presente in

%%%Considerando il problema ottenuto da
%Let us consider the problem obtained from (\ref{P2}) taking
%$\alpha=\alpha_*$ %with $\alpha_*$ in (\ref{2.9})
%\begin{equation}\label{1.3}
%\left\{\begin{array}{l}
%\displaystyle{ v_t \, = \, (a(x)v_x)_x + \alpha_* v \qquad\qquad \, \mbox{in}\,\, Q_T=(-1,1)\times(0,T)
%}\\ [2.5ex]
%\displaystyle{a(x)v_x(t,x)|_{x=\pm 1} = 0 \qquad\qquad\qquad t\in(0,T)}\\ [2.5ex]
%\displaystyle{v (0,x)=v_0(x),\qquad\qquad\,\qquad\qquad\, x\in(-1,1) }~.
%\end{array}\right.
%%\begin{cases}
%%v_t \, = \, (a(x)v_x)_x + \alpha_* v \qquad\qquad \, \mbox{in}\,\, Q_T=(-1,1)\times(0,T)
%%\\
%%a(x)v_x(t,x)|_{x=\pm 1} = 0 \qquad\qquad\qquad t\in(0,T),
%%\\
%%v (0,x)=v_0(x),\qquad\qquad\,\qquad\qquad\, x\in(-1,1)
%%\end{cases}
%\end{equation}
We can see, by Lemma \ref{Autof}, that %an eigenfunction of the equation of (\ref{1.3})
%è la funzione:
%is the function
%$$
%\frac{v_d(x)}{\|v_d\|%_{L^2(\Omega)}
%}\,,
%$$
%%%%%%%%%%%%%%%%%%%%%%%%%%%%%%%%%%%%%%%%%%%%%%%%%%%%%%%%%%%%%%%%%%%%%%%%%%%%%%%%%%%%%%
%%  VERIFICA CHE \alpha_* è un'autofunzione
%%%%%%%%%%%%%%%%%%%%%%%%%%%%%%%%%%%%%%%%%%%%%%%%%%%%%%%%%%%%%%%%%%%%%%%%%%%%%%%%%%%%%%
%Infatti
%In fact, we have
%$$
%\bigg(a(x)\frac{v_{dx}(x)}{\|v_d\|_{L^2(\Omega)}}\bigg)_x +
%\alpha_* \frac{v_d(x)}{\|v_d\|_{L^2(\Omega)}} - \lambda
%\frac{v_d(x)}{\|v_d\|_{L^2(\Omega)}} =
%$$
%$$
%=\frac{(a(x)v_{dx}(x))_x}{\|v_d\|_{L^2(\Omega)}}-\frac{(a(x)v_{dx}(x))_x}{v_d(x)}\frac{v_d(x)}{\|v_d\|_{L^2(\Omega)}}-\lambda
%\frac{v_d(x)}{\|v_d\|_{L^2(\Omega)}}=0\,,
%$$
%da cui
%from which it follows that $\frac{v_d(x)}{\|v_d\|_{L^2(\Omega)}}$
%%è un'autofunzione per
%is an eigenfunction of (\ref{1.3})
%associata all'autovalore
%associated with the eigenvalue $\lambda=0$.
%\\
%%Indichiamo con
%We denote by $\{\omega_k\}^\infty_{k=1}$
%%le autofunzioni ortogonali in
%the orthonormal eigenfunctions in $L^2(\Omega)$
%%del problema
%of the problem (\ref{1.3}),
%%quindi
%Per quanto detto e considerato in particolare che
%%%%%%%%%%%%%%%%%%%%%%%%%%%%%%%%%%%%%%%%%%%%%%%%%%%%%%%%%%%%%%%%%%%%%%%%%%%%%%
%% FINE VERIFICA AUTOFUNZIONE
%%%%%%%%%%%%%%%%%%%%%%%%%%%%%%%%%%%%%%%%%%%%%%%%%%%%%%%%%%%%%%%%%%%%%%%%%%%%%%

\begin{equation}
\label{4}
%\,\exists\, k_* \in \N\,:\qquad\qquad
\lambda_1=0\qquad \mbox{ and } \qquad\omega_{1}(x)=\frac{v_d(x)}{\|v_d\|%_{L^2(\Omega)}
}>0,\,\,\forall
x\in (-1,1)\,.
\end{equation}

\noindent
\underline{STEP.3}
Let us now choose the following static bilinear control
$$
\alpha(x)=\alpha_*(x)+\beta,\,\forall x\in(-1,1),\,\,\mbox{ with } \beta\in \R \mbox{  ($\beta$ to be determined below).}
$$
The corresponding solution of (\ref{P2}),
%per questo particolare
for this particular bilinear coefficient $\alpha,$
%è data da:
has the following Fourier series representation
(\footnote{ Observe that adding $\beta\in\R$ in the coefficient $\alpha_*$ there is a shift of the eigenvalues corresponding to $\alpha_*$ from $\{-\lambda_k\}_{k\in\N}$ to $\{-\lambda_k+\beta\}_{k\in\N},$ but the eigenfunctions remain the same for $\alpha_*$ and $\alpha_*+\beta$. })
$$
v(t,x)=\sum^\infty_{k=1}
e^{(-\lambda_k+\beta)t}\bigg(\int^1_{-1}v_0(s)\omega_k(s)ds\bigg)\omega_k(x)=$$
$$=e^{\beta t}\bigg(\int^1_{-1}v_0(s)\omega_1(s)ds\bigg)\omega_1(x)+\sum_{k>1}
e^{(-\lambda_k+\beta)t}\bigg(\int^1_{-1}v_0(s)\omega_k(s)ds\bigg)\omega_k(x)
$$
%Poniamo
Let % us suppose that
$$
r(t,x)=\sum_{k>1}
e^{(-\lambda_k+\beta)t}\bigg(\int^1_{-1}v_0(s)\omega_k(s)ds\bigg)\omega_k(x)
$$
%dove, tenendo conto che
where, recalling that $\lambda_k%=k(k+1)
<\lambda_{k+1},$
%si ha:
we obtain
$$
-\lambda_k<-\lambda_1=0\quad\mbox{  for ever   } \, k\in \N, \,k>1\,.
$$
%Adesso diamo la seguente stima, tenendo conto della
%Now we give the following estimate,
Owing to (\ref{4}),
$$\|v(t,\cdot)-v_d\|\leq\bigg\|e^{\beta
t}\bigg(\int^1_{-1} v_0(s)\omega_1(s)ds\bigg)\omega_1
-\|v_d\|%_{L^2(\Omega)}
\omega_1\bigg\|%_{L^2(\Omega)}
\!\!\!+\|r(t,x)\|%_{L^2(\Omega)}
\!\!=$$%\\
$$=\left|e^{\beta t}\bigg(\int^1_{-1}
v_0(x)\omega_1(x)dx\bigg)-\|v_d\|%_{L^2(\Omega)}
\right|+\|r(t,x)\|%_{L^2(\Omega)}
$$
%\end{multline}
%Adesso scegliamo
%Since $v_0\in L^2(-1,1), \,v_0\geq 0 \mbox{ and } v_0\not\equiv0$ in $(-1,1)$ and by (\ref{4}),
%%si ha
%we obtain
%\begin{equation}
%\label{6}
%\int^1_{-1} v_0(x)\omega_{1}(x)dx>0.
%\end{equation}
%Then, it is possible choose $\beta$ and $T>0$
%%in modo che
%so that
%$$
%e^{\beta T}\int^1_{-1} v_0\omega_1 dx=\|v_d\|%_{L^2(\Omega)}
%,
%$$
%%cioè
%that is,
%\begin{equation}
%\label{9} \beta =\frac{1}{T}\ln\bigg(
%\frac{\|v_d\|%_{L^2(\Omega)}
%}{\int^1_{-1} v_0\omega_1 dx}\bigg).
%\end{equation}
%%Inoltre possiamo dare la seguente stima, tenendo conto che
Since $-\lambda_k<-\lambda_2,$ $\forall k>2$,
%ed applicando la disuguaglianza di Cauchy-Schwarz
 applying  Parseval's equality we have %can also give the following estimate
%\begin{multline}\label{10}
$$\|r(t,x)\|^2%_{L^2(\Omega)}
\leq e^{2(-\lambda_2 +\beta)t}
\sum_{k>1}\bigg|\int^1_{-1} v_0\omega_k ds\bigg|^2 \|\omega_k(x)\|^2%_{L^2(\Omega)}
=$$
%%\\
$$
=e^{2(-\lambda_2 +\beta)t}\sum_{k>1}\langle v_0,\omega_k\rangle^2= e^{2(-\lambda_2
+\beta)t}\|v_0\|^2%_{L^2(\Omega)}
.$$
%\end{multline}
%Allora da
%So, by (\ref{9}),  the last inequality, and the above estimate for $\|v(T,\cdot)-v_d(\cdot)\|$ we conclude that
%%segue
%%it follows
%$$
%\|v(T,\cdot)-v_d(\cdot)\|\leq e^{(-\lambda_2
%+\beta)T}\|v_0\|%_{L^2(\Omega)}
%%=
%%$$
%%$$
%=e^{-\lambda_2 T}%\ln\bigg(
%\frac{\|v_d\|%_{L^2(\Omega)}
%}{\int^1_{-1} v_0\omega_1
%dx}%\bigg)
%\|v_0\|%_{L^2(\Omega)}
%\stackrel{T\longrightarrow\infty}{\longrightarrow}0\,.
%$$
%%Da cui la tesi, in quanto
%From which %the thesis, being
%we have the conclusion.
%$$
%\forall\, \varepsilon>0\,\exists\, T=T(\varepsilon,v_0,v_d)\mbox{
%and }\alpha=\alpha_*+\beta
%$$
%%con
%with $\alpha_*,$
%%come nella
%as in (\ref{2.9}) and $\beta$
%%come nella
%as in (\ref{9}),
%%tale che
%%so that
%then
%$$
%\|v(T,\cdot)-v_d\|_{L^2(-1,1)}\leq\varepsilon.
%%\qquad\qed
%$$
%This is the conclusion.

%%%%%%%%%%%%%%%%%%%%%%%%%%%%%%%%%%%%%%%%%%%%%%%%%%%%%%%%%%%%%%%%%%%%%%%%%%%%%%%%%%%%%%%%%%%%%%%%%%%%%%%%%%%%%
\noindent Fixed $\varepsilon>0$, we choose $T_\varepsilon >0$ such that
\begin{equation}
\label{T}
e^{-\lambda_2 T_\varepsilon}=\varepsilon\frac{\int^1_{-1} v_0 v_d dx}{\|v_0\|\|v_d\|^2}\,.
\end{equation}
Since $v_0\in L^2(-1,1), \,v_0\geq 0 \mbox{ and } v_0\not\equiv0$ in $(-1,1)$ and by (\ref{4}),
%si ha
we obtain
\begin{equation}
\label{6}
\langle v_0,\omega_{1}\rangle=\int^1_{-1} v_0(x)\omega_{1}(x)dx>0.
\end{equation}
\noindent Then, it is possible choose $\beta_\varepsilon$ %and $T>0$
%in modo che
so that
$$
e^{\beta_\varepsilon T_\varepsilon} \int^1_{-1} v_0\omega_1 dx=\|v_d\|%_{L^2(\Omega)}
\,,
$$
%cioè
that is, since $\omega_1=\frac{v_d}{\|v_d\|},$
\begin{equation}
\label{9} \beta_\varepsilon =\frac{1}{T_\varepsilon}\ln\bigg(
\frac{\|v_d\|^2%_{L^2(\Omega)}
}{%\int^1_{-1} v_0\omega_1 dx
\int^1_{-1} v_0 v_d dx}\bigg).
\end{equation}
So, by (\ref{T}), (\ref{9}) and the above estimates for $\|v(T_\varepsilon,\cdot)-v_d(\cdot)\|$ and $\|r(T_\varepsilon,\cdot)\|$ we conclude that

%So, by $(\ref{v-vd})-(\ref{T})$ %(\ref{9})
%and (\ref{9}) % the last inequality, and the above estimate for $\|v(T,\cdot)-v_d(\cdot)\|$
% we conclude that
%%segue
%%it follows
$$
\|v(T_\varepsilon,\cdot)-v_d(\cdot)\|\leq e^{(-\lambda_2
+\beta_\varepsilon)T_\varepsilon}\|v_0\|%_{L^2(\Omega)}
%=
%$$
%$$
=e^{-\lambda_2 T_\varepsilon}%\ln\bigg(
\frac{\|v_d\|^2%_{L^2(\Omega)}
}{%\int^1_{-1} v_0\omega_1 dx
\int^1_{-1} v_0 v_d dx}%\bigg)
\|v_0\|%_{L^2(\Omega)}
=\varepsilon\,.
%\stackrel{T\longrightarrow\infty}{\longrightarrow}0\,.
$$
%Da cui la tesi, in quanto
From which %the thesis, being
we have the conclusion.
%%%%%%%%%%%%%%%%%%%%%%%%%%%%%%%%%%%%%%%%%%%%%%%%%%%%%%%%%%%%%%%%%%%%%%%%%%%%%%%%%%%%%
\hspace{5.5cm}\eproof }

\vspace{0.5cm}

\proof{(of Theorem \ref{C1})}{
%La dimostrazione del
The proof of Theorem \ref{T1}
%si può adattare, tenendo conto che la precedente
can be adapted to Theorem \ref{C1}, keeping in mind that, in STEP.3, %$\omega_{k_*}(x)=\frac{v_d(x)}{\|v_d(x)\|} > 0$ and
inequality (\ref{6})
%continua a valere in quanto:
continues to hold in this new setting. %also in the Theorem \ref{C1}.
In fact we have
$$
\int^1_{-1} v_0(x)\omega_1 (x) dx=\int^1_{-1}
v_0(x)\frac{v_d(x)}{\|v_d\|}dx=
$$
$$
=\frac{1}{\|v_d\|}\int^1_{-1} v_0v_d dx>0, \mbox{  by assumptions
(\ref{H2}).}
$$
From this point on, %(\ref{7}) continues to hold also in the Theorem \ref{C1},
one can proceed as in the proof of  %STEP.3 and in STEP.4
%della precedente dimostrazione del
%of the previous proof of
Theorem \ref{T1}.
%%Però in queste ipotesi non vale il
%But in this hypothesis, Lemma \ref{NN}
%%quindi la soluzione può anche assumere valori negativi in
%is not valid, so the solution can also assume negative values in $Q_T.$%\qquad\qed
\hspace{3.0cm}\eproof}
%\vspace{1cm}

%\vspace{1.5cm}

%%%%%%%%%%%%%%%%%%%%%%%%%%%%%%%%%%%%%%%%%%%%%%%%%%%%%%%%%%%%%%%%%%%%%%%%%%%%%%%%%

\begin{center} \bf Acknowledgments\end{center}
This research has been performed in the framework of the GDRE CONEDP. The authors wish to thank Institut Henri Poincar\'e (Paris, France) for providing a very stimulating environment during the "Control of Partial and Differential Equations and Applications" program in the Fall 2010.

%\vfill\eject

%\label{lastpage}

\end{document}